\documentclass[12pt,oneside,reqno]{amsart}
\usepackage{amsmath,amssymb,amsthm,textcomp}
\usepackage{amsfonts,graphicx}
\usepackage[mathscr]{eucal}
\usepackage{color}
\usepackage[table]{xcolor}
\usepackage{booktabs}

\usepackage{subfig}
\usepackage{epsfig}
\usepackage{caption}
\theoremstyle{definition}
\usepackage{float}
\usepackage{mathtools}  
\usepackage{diffcoeff}  
\floatstyle{plaintop}
\numberwithin{equation}{section}

\newcommand{\ncom}{\newcommand}

\ncom{\beq}{\begin{equation}}
	\ncom{\eeq}{\end{equation}}
\ncom{\bea}{\begin{eqnarray*}}
	\ncom{\eea}{\end{eqnarray*}}
\ncom{\beqa}{\begin{eqnarray}}
	\ncom{\eeqa}{\end{eqnarray}}
\ncom{\nno}{\nonumber}
\ncom{\non}{\nonumber}
\ncom{\ds}{\displaystyle}
\ncom{\half}{\frac{1}{2}}
\ncom{\mbx}{\makebox{.25cm}}
\ncom{\hs}{\mbox{\hspace{.25cm}}}
\ncom{\rar}{\rightarrow}
\ncom{\Rar}{\Rightarrow}
\ncom{\noin}{\noindent}
\ncom{\bc}{\begin{center}}
	\ncom{\ec}{\end{center}}
\ncom{\sz}{\scriptsize}
\ncom{\rf}{\ref}
\ncom{\s}{\sqrt{2}}
\ncom{\sgm}{\sigma}
\ncom{\Sgm}{\Sigma}
\ncom{\psgm}{\sigma^{\prime}}
\ncom{\dt}{\delta}
\ncom{\Dt}{\Delta}
\ncom{\lmd}{\lambda}
\ncom{\Lmd}{\Lambda}
\ncom{\Th}{\Theta}
\ncom{\e}{\eta}
\ncom{\eps}{\epsilon}
\ncom{\pcc}{\stackrel{P}{>}}
\ncom{\lp}{\stackrel{L_{p}}{>}}
\ncom{\dist}{{\rm\,dist}}
\ncom{\sspan}{{\rm\,span}}
\ncom{\re}{{\rm Re\,}}
\ncom{\im}{{\rm Im\,}}
\ncom{\sgn}{{\rm sgn\,}}
\ncom{\ba}{\begin{array}}
	\ncom{\ea}{\end{array}}
\ncom{\hone}{\mbox{\hspace{1em}}}
\ncom{\htwo}{\mbox{\hspace{2em}}}
\ncom{\hthree}{\mbox{\hspace{3em}}}
\ncom{\hfour}{\mbox{\hspace{4em}}}
\ncom{\vone}{\vskip 2ex}
\ncom{\vtwo}{\vskip 4ex}
\ncom{\vonee}{\vskip 1.5ex}
\ncom{\vthree}{\vskip 6ex}
\ncom{\vfour}{\vspace*{8ex}}
\ncom{\norm}{\|\;\;\|}
\ncom{\integ}[4]{\int_{#1}^{#2}\,{#3}\,d{#4}}
\ncom{\vspan}[1]{{{\rm\,span}\{ #1 \}}}
\ncom{\dm}[1]{ {\displaystyle{#1} } }
\ncom{\ri}[1]{{#1} \index{#1}}

\newtheorem{theorem}{\bf Theorem}[section]
\newtheorem{remark}{\bf Remark}[section]
\newtheorem{proposition}{Proposition}[section]
\newtheorem{lemma}{Lemma}[section]
\newtheorem{corollary}{Corollary}[section]
\newtheorem{example}{Example}[section]

\newtheoremstyle
{remarkstyle}
{}
{11pt}
{}
{}
{\bfseries}
{:}
{     }
{\thmname{#1} \thmnumber{#2} }

\theoremstyle{remarkstyle}

\def\eps{\varepsilon}

\begin{document}
	\title{\large A p\lowercase{robabilistic}  E\lowercase{xtension}  \lowercase{of~the} f\lowercase{ubini} p\lowercase{olynomials}} 
		\author[Ritik Soni]{R. Soni$^{1}$}
	\author{A. K. Pathak$^{1*}$}
	\author{P. Vellaisamy$^{2}$\\
		$^{1}$D\lowercase{epartment of} M\lowercase{athematics and} S\lowercase{tatistics}, C\lowercase{entral} U\lowercase{niversity of} P\lowercase{unjab},\\
		B\lowercase{athinda}, P\lowercase{unjab}-151401, I\lowercase{ndia}.
		\\$^{2}$D\lowercase{epartment of} M\lowercase{athematics}, I\lowercase{ndian} I\lowercase{nstitute of} T\lowercase{echnology}
		B\lowercase{ombay},\\ P\lowercase{owai}, M\lowercase{umbai}-400076, I\lowercase{ndia}.}
\thanks{*Corresponding Author}
\thanks{ E-mail Address: ritiksoni2012@gmail.com (R. Soni), ashokiitb09@gmail.com (A. K. Pathak),} \thanks{pv@math.iitb.ac.in (P. Vellaisamy)}
\thanks{The research of  R. Soni was supported by CSIR, Government of India.}

\begin{abstract} In this paper, we  present a probabilistic extension of the Fubini polynomials and numbers associated with a random variable  satisfying some appropriate moment conditions. We obtain the exponential generating function and an integral representation for it. The higher order Fubini polynomials and recurrence relations are also derived.  A probabilistic generalization of a series transformation formula and some interesting examples are discussed.  A connection between the probabilistic Fubini polynomials and  Bernoulli, Poisson, and geometric random variables are also established.  Finally, a determinant expression formula is presented.

\end{abstract}

\maketitle
\section{Introduction} \noindent Recently, certain polynomials and numbers have received growing attention in many branches of mathematics, computer science, and physics. More specifically, the study on the Stirling numbers of the second kind has progressed significantly during the past few decades. The Stirling numbers of the second kind, denoted by $S(n,k)$, count the total number of partitions of a set of $n$ elements into $k$ non-empty disjoint subsets and  play an important role in combinatorics. It is defined by 
(see \cite{quaintance2015combinatorial})
\begin{equation}\label{0004}
S(n,k) = \frac{1}{k!} \sum_{j=0}^{k} (-1)^{k-j}\binom{k}{j}j^n.
\end{equation}
Its exponential generating function is given by 
(see \cite[Chapter 9]{quaintance2015combinatorial})
\begin{equation}\label{0003}
\sum_{n=k}^{\infty} S(n,k) \frac{t^n}{n!} = \frac{(e^t-1)^k}{k!} ,\;\; t \in \mathbb{C},
\end{equation}
where $\mathbb{C}$ is the set of complex numbers.
For more details on the Stirling numbers of the second kind and its properties, one may refer to Comtet \cite{comtet2012advanced} and Gould \cite{quaintance2015combinatorial}.

\vspace*{.3cm}
\noindent A number of polynomials are defined through $S(n,k)$.  For instance, the Bell polynomials $B_n(x)$ are defined as  (see \cite{kataria2022probabilistic} and \cite{kim2019central})
\begin{equation*}
B_n(x) = \sum_{k=0}^{n}S(n,k)x^k.
\end{equation*}
An alternate expression for  $B_n(x)$ is 
\begin{equation}\label{00001}
B_n(x) = \sum_{k=0}^{\infty}k^n  \frac{e^{-x}x^k}{k!}= E(Y^n(x)),
\end{equation} 
which is the $n$-th moment of the Poisson variable $Y(x)$ with mean $x>0$.

\noindent Besides the Bell polynomials, the Fubini polynomials are also applied in various disciplines of the applied sciences and   combinatorics. These polynomials are also known as the geometric polynomials or the ordered Bell polynomials. Through the  Stirling numbers of the second kind,  Fubini polynomials are defined by the relation (see \cite{adell2022higher,  boyadzhiev2016geometric, qi2019determinantal, tanny1975some})
\begin{equation}\label{112}
W_n(x) = \sum_{k=0}^{n} k! S(n,k)x^k.
\end{equation}
The exponential generating function of $W_n(x)  $  is
\begin{equation}\label{gf1}
 \sum_{n=0}^{\infty} W_n(x) \frac{t^n}{n!} = \Big[1-x\left(e^{t}-1 \right) \Big]^{-1}. 
\end{equation}
Note that, when $x=1$, (\ref{112}) yields 
\begin{equation}\label{111}
W_n = W_n(1)= \sum_{k=0}^{n} S(n,k)k!,
\end{equation}
which are called  Fubini numbers 
(see \cite{belbachir2021central, boyadzhiev2005series}) and satisfy the  recurrence relation (see \cite{gross1962preferential})
\begin{equation}\label{rr12}
W_n = \sum_{k=1}^{n} \binom{n}{k} W_{n-j}.
\end{equation}

\vspace*{.3cm}
Recently, Adell and Lekuona \cite{adell2019probabilistic} defined a probabilistic version of the Stirling numbers of the second kind. Let $Y$ be a real valued random variable (rv) having finite moment generating function ($mgf$) and  $\{ Y_j \}_{j \geq 0}$ be a sequence of independent and identically distributed (IID) random variables (rvs) with distribution as that of the rv $Y.$ Define $S_0 =0$ and $S_j = Y_1 +Y_2 + \cdots + Y_j, ~j \geq 1.$  The probabilistic  Stirling numbers of the second kind $S_Y(n,k)$, associated with the  rv $Y$, is defined via the relation
\begin{equation}\label{0002}
S_Y(n,k) = \frac{1}{k!} \sum_{j=0}^{k} (-1)^{k-j}\binom{k}{j}\mathbb{E}S_j^n.
\end{equation}
Its exponential generating function is (see \cite{adell2020probabilistic} and \cite{adell2019probabilistic})
\begin{equation}\label{0001}
 \sum_{n=k}^{\infty} S_Y(n,k) \frac{t^n}{n!}= \frac{\left(\mathbb{E}e^{tY}-1\right)^k}{k!},\;\;t \in \mathbb{C}.
\end{equation}
When $Y$ is degenerate at 1,  (\ref{0002}) and (\ref{0001}) reduces to (\ref{0004}) and (\ref{0003}), respectively. They obtained the moments of 
$S_j$ as (see \cite{adell2020probabilistic}) 
\begin{equation}\label{def1a}
\mathbb{E}S_j^n = \sum_{k=0}^{n \wedge j} S_Y(n,k) (j)_k,
\end{equation}
where $n \wedge j = \min\{n,j\}$ and $(j)_k = j(j-1)\cdots(j-k+1)$ is the falling factorial.\\
Soni {\it et al.} \cite{rsoni2022pbp} discussed the probabilistic Bell polynomials defined by
\begin{equation}\label{bp11}
B_n^Y(x) = \sum_{k=0}^{n}S_Y(n,k)x^k,
\end{equation}
which has the exponential generating function
\begin{equation}\label{gfpbp1}
e^{\left(\mathbb{E}e^{tY}-1\right)} = \sum_{n=0}^{\infty} B_n^Y(x) \frac{t^n}{n!}.
\end{equation}
An alternative representation of $B_n^Y(x)$ in terms of the Poisson moments is given by
\begin{equation}\label{0009}
B_n^Y(x) = \sum_{k=0}^{\infty}\mathbb{E}S_k^n \frac{e^{-x}x^k}{k!}.
\end{equation}
When $Y$ is degenerate at $1$, it coincides with (\ref{00001}),  the famous Dobi\'{n}ski's formula.

\vspace*{.3cm}
\noindent Kim \cite{kim2022degenerate} established a connection between the geometric rv and the Fubini polynomials.  For  $p \in (0,1]$,  let  $G_p$ be a geometric rv with probability mass function ($pmf$) $P\{G_p=i\}=p(1-p)^{i-1},\;i \ge 1.$ Note that  $G_p$  denotes the number of trials needed for the first success, when a coin with success probability $p$ is tossed. It follows easily that
\begin{align}\label{grv1}
\mathbb{E}[G_p-1]^n =& \sum_{i=1}^{\infty} (i-1)^n P\{G_p=i\} \nonumber \\
 =& p\sum_{i=1}^{\infty} (i-1)^n (1-p)^{i-1} \nonumber \\
  =& p\sum_{i=0}^{\infty} i^n (1-p)^{i}.
\end{align}
For $x\geq 0$, let $\eta(x)=(1+x)^{-1}.$ Then the connection between the geometric rv $G_{\eta (x)}$ and the Fubini polynomials is
\begin{align}\label{gg}
W_n(x) =& \sum_{i=0}^{\infty}i^n P\{G_{\eta_{(x)}} =i+1\} 
         =\mathbb{E}[G_{\eta_{(x)}}-1]^n. 
\end{align}

\vspace*{.3cm}
As mentioned earlier, a probabilistic representation of the  Stirling numbers of the second kind in terms of IID rvs is  studied by Adell and Lekuona \cite{adell2019probabilistic} and Adell \cite{ adell2020probabilistic}. These results are very useful in the analytical number theory and  in generalizing different classical sums of powers of arithmetic progression formulas. Laskin \cite{laskin2009some} studied  a fractional generalization of the Bell polynomials and the Stirling numbers of the second kind.  Guo and Zhu \cite{guo2020generalized} introduced the generalized Fubini polynomials and studied their logarithmic properties.
Motivated by the work of Adell and Lekuona \cite{adell2019probabilistic} and 
  Guo and Zhu \cite{guo2020generalized}, we consider a probabilistic generalization of the Fubini polynomials and numbers and explore their important properties. The connections of these polynomials with the known families of the probability distributions are explored. A simple determinant formula for the probabilistic Fubini numbers is also derived along with some combinatorial sums.

\vspace*{.3cm}
The paper is organized as follows. In Section 2, we present a probabilistic extension of the Fubini polynomials and numbers. The exponential generating function, integral representation, and some recurrence relations are obtained. A connection between the higher order probabilistic Fubini polynomials and the negative binomial process is also discussed. In Section 3, we obtain a probabilistic generalization of a series transformation formula and illustrate it with some examples. A new relationship between the rising factorial and the Lah numbers is deduced. A connection of the probabilistic Fubini polynomials with  Bernoulli,  Poisson, and geometric random variates are discussed in Section 4. Finally, a determinant expression for the probabilistic Fubini numbers and  a combinatorial sum formula is also obtained in Section 5.

\section{Probabilistic Fubini Polynomials and Numbers }
\noindent Let $\mathcal{G}$ be the set of rvs $Y$ satisfying the following moment conditions
\begin{equation}\label{cond1}
\mathbb{E}|Y|^n < \infty, \;\;\;\; n \in \mathbb{N}_0, \;\;\;\; \lim_{n \rightarrow \infty} \frac{|t|^n \mathbb{E}|Y|^n}{n!} =0,\;\;\; |t| <r,
\end{equation}
where $\mathbb{N}_0 = \mathbb{N} \cup \{0\}$, $r > 0,$ $\mathbb{N}$ is the set of natural numbers and $\mathbb{E}$ denotes the mathematical expectation.
Equation (\ref{cond1}) confirms the existence of the moment generating function for rv $Y$ (see \cite[p. 344]{billingsley2008probability}).

Let $\langle Y_i\rangle_{i \geq 0}$ be  independent and identically distributed (IID) copies of a rv $Y \in \mathcal{G}.$
By Jensen's inequality, we have 
\begin{equation*}
\mathbb{E}|S_i|^n \leq i^{1/n} \mathbb{E}|Y|^n, ~~\text{ for }  ~~n \in \mathbb{N},
\end{equation*} 
where $S_i = Y_1 + Y_2 + \cdots + Y_i$  and $S_0 =0$.\\

In view of \eqref{gg}, we define the probabilistic Fubini polynomials associated with the rv $Y$ as
\begin{equation}\label{ad}
W_n^Y(x) = \sum_{i=0}^{\infty} \mathbb{E}\left[S_i^n\right] P\{G_{\eta(x)} =i+1\},
\end{equation}
where $G_p$ follows the geometric distribution with parameter $p$ and $\eta(x)=(1+x)^{-1}.$
In case of the degeneracy of the rv $Y$ at 1, (\ref{ad}) leads to the classical Fubini polynomials. In particular, when $Y=1$ and $x=1$ in (\ref{ad}), we get the $n$th order moments of the geometric rvs with probability of success equals 1/2. These moments are well-known as the Fubini numbers.

\vspace*{.3cm}
Next, we present the exponential generating function for the probabilistic Fubini polynomials.
\begin{proposition}\label{p2}
Let $|x\left(\mathbb{E}e^{tY}-1\right)| \leq 1$. Then,
\begin{equation}\label{pp1}
\frac{1}{1-x\left(\mathbb{E}e^{tY}-1\right)} = \sum_{n=0}^{\infty} W_n^Y(x) \frac{t^n}{n!}.
\end{equation}
\end{proposition}
\begin{proof}
For the IID copies of the rv $Y$,  we have, from (\ref{ad}),
\begin{align*}
\sum_{n=0}^{\infty} W_n^Y(x) \frac{t^n}{n!} &=  \sum_{n=0}^{\infty} \left(\sum_{k=0}^{\infty} \mathbb{E}S_k^n \left(\frac{1}{1+x}\right) \left(\frac{x}{1+x}\right)^k\right) \frac{t^n}{n!}\\
&= \sum_{k=0}^{\infty} \left(\frac{1}{1+x}\right) \left(\frac{x}{1+x}\right)^k \left(\mathbb{E}e^{tS_k}\right)\\
&= \sum_{k=0}^{\infty} \left(\frac{1}{1+x}\right) \left(\frac{x}{1+x}\right)^k \left(\mathbb{E}e^{tY}\right)^k\\
&=	\frac{1}{1-x\left(\mathbb{E}e^{tY}-1\right)}.
\end{align*}
Hence, the proposition is proved.
\end{proof}
The series expansion of (\ref{pp1}) in the light of (\ref{0001}) gives an alternative  representation of the  probabilistic Fubini polynomials in terms of the  probabilistic Stirling numbers of the second kind of the following form 
\begin{equation}\label{113}
W_n^Y(x) = \sum_{k=0}^{n} S_Y(n,k)k!x^k.
\end{equation}

It may be observed that for $x=1$, (\ref{113}) gives a probabilistic generalization of the classical Fubini numbers given by
\begin{equation}\label{fn131}
W_n^Y = \sum_{k=0}^{n} S_Y(n,k)k!.
\end{equation}
We call it the probabilistic Fubini numbers.

When $Y$ follows an exponential distribution with mean 1, we establish a connection between the probabilistic 
Fubini polynomials and the Lah numbers of the following form
\begin{equation}\label{lah1}
W_n^Y(x) = \sum_{k=0}^{n} L(n,k) k! x^k.
\end{equation}
In the literature, (\ref{lah1}) is termed as  $0$-Fubini-Lah polynomials (see \cite{racz2020r}).

In the following proposition,  we give an integral representation for the probabilistic Fubini polynomials.
\begin{proposition}\label{p1}
Let $Y \in \mathcal{G}$. Then
\begin{equation}\label{exp12}
W_n^Y(x) = \int_{0}^{\infty} B_n^Y(xv)e^{-v}d v = \mathbb{E}_V[B_n^Y(xV)],
\end{equation}
where $V$ is the standard exponential rv with probability density function $f(v) = e^{-v}, v\geq 0$ and $\mathbb{E}_V$ stands for the mathematical expectation for rv $V$.

\begin{proof}
Considering (\ref{113}) and with the help of the gamma integral, we have
\begin{align*}\label{int1}
\mathbb{E}_V[B_n^Y(xV)]  &=\int_{0}^{\infty} B_n^Y(xv)e^{-v}d v \\
&=  \int_{0}^{\infty} \left( \sum_{k=0}^{n}S_Y(n,k)x^k v^k \right)e^{-v}d v\\
&= \sum_{k=0}^{n}S_Y(n,k)k! x^k \\
&= W_n^Y(x),
\end{align*}
where $k! = \int_{0}^{\infty} v^k e^{-v}dv$.
\end{proof}
\end{proposition}
For $n$ and $k$ be two non-negative integers such that $n \geq k,$ the partial exponential Bell polynomials  $B_{n,k}(x_1,x_2,\dots, x_{n-k+1})$ have the following form  (see \cite{comtet2012advanced})
\begin{equation*}
	B_{n,k}(x_1 , x_2, x_3,...., x_{n-k+1} ) = n!\biggl(\sum_{\Lambda_n ^{k}}^{} \prod_{j=1}^{n-k+1} \frac{1}{k_j !} \biggl(\frac{x_j}{j!}\biggr)^{k_j}\biggr),
\end{equation*}
where the summation is taken over the following set
\begin{equation*}
	\Lambda_n ^{k} = \biggl\{(k_1, k_2, ...., k_{n-k+1}) : \sum_{j=1}^{n-k+1}k_j = k, \sum_{j=1}^{n-k+1}jk_j = n, k_j 	\in\mathbb{N}_0 \biggr\}.
\end{equation*}
A connection between the probabilistic Stirling numbers of the second kind and the partial exponential Bell polynomials is obtained and is given by (see \cite{rsoni2022pbp})
\begin{equation}\label{snpb1}
	S_Y(n,k) =  B_{n,k} \left(\mathbb{E}Y, \mathbb{E}Y^2, \dots,\mathbb{E}Y^{n} \right).
\end{equation}
\noindent When $x=1$, and using (\ref{snpb1}) and  Proposition \ref{p1}, we get
\begin{equation*}
W_n^Y = \mathbb{E}_V  B_{n} \left(V\mathbb{E}Y, V\mathbb{E}Y^2, \dots, V\mathbb{E}Y^{n} \right),
\end{equation*}
where $B_{n}$ are the complete exponential Bell polynomials which can be expressed in terms of the partial exponential Bell polynomials as (see \cite[p. 133]{comtet2012advanced})
\begin{equation*}
B_{n}(x_1,x_2,\dots, x_n) = \sum_{k=0}^{n} B_{n,k}(x_1,x_2,\dots, x_{n-k+1}).
\end{equation*}

It is well-known that geometric distribution is a special case of the negative binomial distribution. For $\alpha > 0$ and $0< p < 1$, let $Z_p$ follows negative binomial distribution denoted by NB$(\alpha, p)$ with $pmf$ 
\begin{equation*}
P\{Z_p = i \} = \binom{-\alpha}{i} \left(p -1 \right)^i p^{\alpha},\;\;i \in \mathbb{N}_0.
\end{equation*}
When $ p =\eta(x)= 1/(1+x)$, the  $mgf$ of $Z_{\eta(x)}$ is given by
\begin{equation*}
\mathbb{E}e^{tZ_{\eta(x)}} = \sum_{i=0}^{\infty} \binom{-\alpha}{i} \left(-\frac{e^tx}{1+x}\right)^i \left(\frac{1}{1+x}\right)^\alpha = \frac{1}{\left(1-x\left(e^t -1\right)\right)^\alpha},
\end{equation*}
provided $t < \log \left(1+1/x\right)$.\\
We define the $\alpha$-th order probabilistic Fubini polynomials as
\begin{equation*}
W_n^{Y}(x;\alpha) = \sum_{i=0}^{\infty} \mathbb{E}\left[S_i^n\right] P\{Z_{\eta(x)} =i\}, \;\;\;\alpha \in \mathbb{N}_0.
\end{equation*}
It has following exponential generating function
\begin{equation}\label{gfho1}
  \sum_{n=0}^{\infty} W_n^{Y}(x;\alpha) \frac{t^n}{n!}=\frac{1}{(1-x\left(\mathbb{E}e^{tY}-1\right))^{\alpha}}.
\end{equation}
Using the series expansion formula $\frac{1}{(1-x)^\alpha} = \sum_{i = 0}^{\infty} \binom{-\alpha}{i} (-x)^i$, the exponential generating function (\ref{gfho1}) is simplified as
\begin{align*}
\frac{1}{(1-x\left(\mathbb{E}e^{tY}-1\right))^{\alpha}} &= \sum_{i = 0}^{\infty}  (-x)^i\binom{-\alpha}{i} \left(\mathbb{E}e^{tY}-1\right)^i\\
&= \sum_{i = 0}^{\infty}  (-x)^i\binom{-\alpha}{i} \sum_{n =i}^{\infty} S_Y(n,i) \frac{t^n}{n!}, \;\;\; (\text{using }  (\ref{0001}))\\
&=\sum_{n =0}^{\infty} \sum_{i = 0}^{n} \binom{\alpha + i-1}{i} i! x^i S_Y(n,i) \frac{t^n}{n!}.
\end{align*} 
On comparing with (\ref{gfho1}), we get
\begin{equation}\label{ar11}
W_n^{Y}(x;\alpha) = \sum_{i = 0}^{n} \binom{\alpha + i-1}{i} i! x^i S_Y(n,i),
\end{equation}
which can be viewed as an alternate representation for the $\alpha$th order probabilistic Fubini polynomials. We  also obtain some identities and interconnections of $\alpha$th order probabilistic Fubini polynomials in the subsequent sections.

Next, we obtain some recurrence relations for the probabilistic Fubini polynomials and also discuss their special cases.
\begin{proposition}\label{p3} Let $W_n^{Y}(x)$ be the probabilistic Fubini polynomials. Then, we have 
\begin{equation}\label{rec1}
W_n^{Y}(x) = x\sum_{k=1}^{n} \binom{n}{k} \mathbb{E}Y^k W_{n-k}^Y(x).
\end{equation}	
\end{proposition}
\begin{proof}
Using (\ref{pp1}), we get
\begin{align*}
\sum_{n=1}^{\infty} W_n^Y(x) \frac{t^n}{n!} 
&= \frac{1}{1-x\left(\mathbb{E}e^{tY}-1\right)} -1 = \frac{x\left(\mathbb{E}e^{tY}-1\right)}{1-x\left(\mathbb{E}e^{tY}-1\right)}\\
&= x\left(\sum_{k=0}^{\infty} \mathbb{E}Y^k \frac{t^k}{k}\right)\left( \sum_{n=k}^{\infty} W_{n-k}^Y(x)\frac{t^{n-k}}{(n-k)!}\right) - x\left(\sum_{n=0}^{\infty} W_n^Y(x) \frac{t^n}{n!}\right)\\
&= \sum_{n=1}^{\infty}\left(x\sum_{k=1}^{n} \binom{n}{k} \mathbb{E}Y^k W_{n-k}^Y(x)\right)\frac{t^n}{n!}.
\end{align*}
Comparing the coefficients of $t^n$ on both sides, we get required result.
\end{proof}
With the help of the Proposition \ref{p3}, one can deduce recurrence relation for the probabilistic Fubini numbers. In particular,  when $Y=1$, it reduces to (\ref{rr12}).

\begin{theorem}\label{thm9}
Let $Y \in \mathcal{G}$. Then, for $n \geq k \geq i$, we have
\begin{equation*}
W_{n+1}^Y(x) = x \sum_{k=0}^{n}  \binom{n}{k}  \mathbb{E}Y^{n-k+1} \sum_{i=0}^{k} \binom{k}{i}  W_i^Y(x) W_{k-i}^Y(x).
\end{equation*}
\end{theorem}
\begin{proof}
Differentiating (\ref{pp1}) with respect to $t$ on both sides, we get
\begin{align*}
\sum_{n=0}^{\infty} W_{n+1}^Y(x) \frac{t^n}{n!} 
&= \frac{d}{dt}\left[\frac{1}{1-x\left(\mathbb{E}e^{tY}-1\right)}\right]\\
&= x\frac{\mathbb{E}\left[Ye^{tY}\right]}{\left[1-x\left(\mathbb{E}e^{tY}-1\right)\right]^2}\\
&= x \left(\sum_{i=0}^{\infty} W_{i}^Y(x) \frac{t^i}{i!}\right) \left(\sum_{k=0}^{\infty} W_{k}^Y(x) \frac{t^k}{k!}\right)\left(\sum_{n=0}^{\infty}  \mathbb{E}Y^{n+1} \frac{t^n}{n!}\right)\\
&= \sum_{n=0}^{\infty} \left(x \sum_{k=0}^{n} \binom{n}{k}  \mathbb{E}Y^{n-k+1} \sum_{i=0}^{k} \binom{k}{i}  W_i^Y(x) W_{k-i}^Y(x)\right) \frac{t^n}{n!}.
\end{align*}
Equating the coefficients of $t^n$, we get the result.
\end{proof}
Next, the following proposition is a recurrence relation in terms of derivative for the probabilistic Fubini polynomials. 
\begin{theorem}
For $Y \in \mathcal{G}$, we have 
\begin{equation} \label{thm22}
\frac{d}{dx } W_n^Y(x) =  \sum_{k=0}^{n} \binom{n}{k}  \mathbb{E}Y^{n-k} \sum_{i=0}^{k} \binom{k}{i}  W_i^Y(x) W_{k-i}^Y(x)  - W_k^Y(x) W_{n-k}^Y(x).
\end{equation}
\end{theorem}
\begin{proof}
Differentiating (\ref{pp1}) with respect to $x$, we get
\begin{align*}
\frac{d}{dx }\sum_{n=0}^{\infty} W_n^Y(x) \frac{t^n}{n!}
&= \frac{\mathbb{E}e^{tY}-1}{(1-x\left(\mathbb{E}e^{tY}-1\right))^2}\\
&= \frac{1}{(1-x\left(\mathbb{E}e^{tY}-1\right))^2} \mathbb{E}e^{tY} -  \frac{1}{(1-x\left(\mathbb{E}e^{tY}-1\right))^2}.
\end{align*}
Making series expansion of $1/(1-x\left(\mathbb{E}e^{tY}-1\right))^2$ and with the help of Theorem \ref{thm9}, we get
\begin{align*}
\frac{d}{dx }\sum_{n=0}^{\infty} W_n^Y(x) \frac{t^n}{n!}
=\frac{1}{(1-x\left(\mathbb{E}e^{tY}-1\right))^2} \mathbb{E}e^{tY} -  \left(\sum_{k=0}^{\infty} W_{k}^Y(x) \frac{t^k}{k!}\right) \left(\sum_{n=0}^{\infty} W_{n}^Y(x) \frac{t^n}{n!}\right)\\
=  \sum_{n=0}^{\infty} \sum_{k=0}^{n} \binom{n}{k}  \mathbb{E}Y^{n-k} \sum_{i=0}^{k} \binom{k}{i}  W_i^Y(x) W_{k-i}^Y(x) \frac{t^n}{n!} - \sum_{n=0}^{\infty}\sum_{k=0}^{n} \binom{n}{k} W_k^Y(x) W_{n-k}^Y(x)\frac{t^n}{n!}\\
= \sum_{n=0}^{\infty} \sum_{k=0}^{n} \binom{n}{k}  \mathbb{E}Y^{n-k} \sum_{i=0}^{k} \binom{k}{i}  W_i^Y(x) W_{k-i}^Y(x)  - W_k^Y(x) W_{n-k}^Y(x)\frac{t^n}{n!}.
\end{align*}
On comparing the coefficients of $t^n$ on both sides, the result in \eqref{thm22} follows. 
\end{proof}
\section{Probabilistic Generalization of a Series Transformation Formula}
Spivey \cite{spivey2007combinatorial} recently unveiled a new approach to evaluate combinatorial sums formula using a finite difference technique. These combinatorial sums can be obtained in terms of the Stirling numbers of the second kind. Adell and Lekuona \cite{adell2017note} and Adell \cite{adell2019probabilistic}, studied the applications of the probabilistic Stirling numbers of the second kind and obtained the probabilistic extension of some known combinatorial identities. Boyadzhiev \cite{boyadzhiev2005series} considered  a series transformation formula with numerous examples.
Let $f(x)$ and $g(x)$ be two arbitrary functions such that $f(x)$ is entire and $g(x)$ is analytic on $ D = \{x, r < |x| < R\}$ with $ 0 \leq r < R$. Then, $f(x)$ and $g(x)$ can be written as
\begin{equation}\label{fg111}
f(x) = \sum_{n=0}^{\infty}f_n x^n, \;\;\;\;\;g(x) = \sum_{n=- \infty}^{\infty}g_n x^n,
\end{equation} 
where $f_n$ and $g_n$ denotes the $n$th coefficient of series for the functions $f$ and $g$, respectively.

Motivated by Boyadzhiev \cite{boyadzhiev2005series} work, we present, in the next result, a probabilistic generalization of the series transformation formula. An important feature of this generalization is that for an appropriate choice of the functions $f$ and $g$, and for a suitable rv $Y$ in the class $\mathcal{G}$, several well-known series sums formulas can be obtained in the closed forms involving some known classical polynomials and probability distribution functions. 
\begin{theorem}
For $n \leq i$, we have
\begin{equation}\label{ide}
\sum_{i=0}^{\infty} g^{(i)}(0)\mathbb{E}f(S_i)\frac{x^i}{i!} = \sum_{n=0}^{\infty} \frac{f^{(n)}(0)}{n!} \sum_{k=0}^{n} S_Y(n,k)g^{(k)}(x)x^k,
\end{equation}
where $S_i$ is sum of IID copies of the rv $Y \in \mathcal{G}.$
\end{theorem}
\begin{proof}
Using (\ref{def1a}) for $n \leq i$ , we get
\begin{align*}
\sum_{i= 0}^{\infty}\frac{ g^i(0)}{i!} x^i \mathbb{E}S_i^n &=  \sum_{i=0}^{\infty}\frac{ g^i(0)}{i!} x^i \left(\sum_{k=0}^{ n}\binom{i}{k}k! S_Y(n,k)\right)\\
&= \sum_{k=0}^{ n} S_Y(n,k)  \sum_{i=0}^{\infty} \binom{i}{k} k! \frac{ g^i(0)}{i!} x^i\\
&= \sum_{k=0}^{ n} S_Y(n,k) x^k \sum_{i= k}^{\infty} \frac{ g^i(0)}{i!} i! \frac{x^{i-k}}{(i-k)!}\\
&= \sum_{k=0}^{ n} S_Y(n,k) x^k g^{(k)}(x).
\end{align*}

On multiplying both sides with $ {f^{(n)}(0)}/{n!}$ and summing over $n$ from $0$ to $\infty$, we get the result.
\end{proof}

\begin{remark}
When $Y $ is degenerate at $1$,  (\ref{ide}) reduces to the following  series transformation formula (see \cite[Eq. 4.11]{boyadzhiev2005series})
\begin{equation}\label{stf11}
\sum_{i=0}^{\infty} g^{(i)}(0)f(i)\frac{x^i}{i!} = \sum_{n=0}^{\infty} \frac{f^{(n)}(0)}{n!} \sum_{k=0}^{n} S(n,k)g^{(k)}(x)x^k.
\end{equation}
\end{remark}
\begin{example}
For $g(x)=e^x$ in (\ref{ide}), we get a new identity which is given by
\begin{equation}\label{98}
\sum_{i=0}^{\infty} \mathbb{E}f(S_i)\frac{x^i}{i!} 
= \sum_{n=0}^{\infty} \frac{f^{(n)}(0)}{n!} \sum_{k=0}^{n} S_Y(n,k)e^x x^k
= e^x \sum_{n=0}^{\infty} \frac{f^{(n)}(0)}{n!} B_n^Y(x).
\end{equation}
On rearrangement of terms in (\ref{98}), we obtain a connection of the probabilistic Bell polynomials with Poisson rv of the form
\begin{equation}\label{99}
\sum_{i=0}^{\infty} \mathbb{E}f(S_i) P\{Y(x) = i\} 
= \sum_{n=0}^{\infty} \frac{f^{(n)}(0)}{n!} B_n^Y(x),
\end{equation}
where $Y(x)$ follows  Poisson distribution with parameter $x$.

Moreover, for $f(x)=x^n$, (\ref{99}) leads to the probabilistic Bell polynomials defined in (\ref{0009}).
\end{example}
\begin{corollary}
If $f(x)$ be the polynomial of degree $n$, then
\begin{equation*}
\sum_{n=0}^{\infty} \frac{f^{(n)}(0)}{n!} B_n(x) = \sum_{k=0}^{n} \frac{x^k}{k!} \sum_{j=0}^{k} (-1)^j \binom{k}{j} f(k-j),
\end{equation*}
where $B_n(x)$ are the Bell polynomials.
\end{corollary}
Proof of the corollary can be executed by the idea of (\ref{99}) and Theorem 9.2 of \cite{quaintance2015combinatorial}.

\begin{example}
For $g(x) = \frac{1}{1-x}$ with $|x| < 1$ and $f(x) = x^n$, from (\ref{ide}) we get
\begin{equation}\label{eq1}
\sum_{i=0}^{\infty} \mathbb{E}S_i^n x^i = \frac{1}{1-x} W_n^Y \left(\frac{x}{1-x}\right) = \sum_{k=0}^{\infty}x^k \sum_{j=0}^{n \wedge k}\binom{k}{j} j! S_Y(n,j),
\end{equation}
where $n \wedge k= \min(n,k)$. This is a probabilistic generalization of the following identity studied in \cite{boyadzhiev2005series}.

\begin{equation*}
\sum_{i=0}^{\infty} i^n x^i = \frac{1}{1-x} W_n \left(\frac{x}{1-x}\right).
\end{equation*}
It may be observed that for  $x=\frac{1}{2}$, (\ref{eq1}) gives 
\begin{equation}\label{ind1}
W_n^Y(1) = \frac{1}{2} \sum_{i=0}^{\infty} \frac{\mathbb{E}S_i^n}{2^i},
\end{equation} 
which is a  probabilistic extension of the Fubini numbers.

Also, when $Y$  follows a standard exponential distribution, we have a new connection between $n$-th sum of rising factorial and the Lah numbers, which is given as
\begin{equation*}
\frac{1}{2} \sum_{i=0}^{\infty} \frac{\langle i \rangle _n}{2^i} = \sum_{k=0}^{n} L(n,k) k!.
\end{equation*}
\end{example}

\begin{example}
For $g(x)= \frac{1}{(1-x)^r}$ with $Re(r) >0$ and $|x| < 1,$ (\ref{ide}) yields 
\begin{equation}\label{ni1}
\sum_{i=0}^{\infty} (-1)^i\binom{-r}{i} \mathbb{E}f(S_i) x^i = \sum_{n=0}^{\infty}  \frac{f^{(n)}(0)}{n!} \sum_{k=0}^{n} S_Y(n,k) (k+r-1)!  x^k \frac{1}{(1-x)^{k+r}}.
\end{equation}
Alternatively, (\ref{ni1}) may be expressed as
\begin{equation}\label{ni11}
\sum_{i=0}^{\infty} (-1)^i\binom{-r}{i} \mathbb{E}f(S_i) x^i = \frac{1}{(1-x)^{r}}\sum_{n=0}^{\infty}  \frac{f^{(n)}(0)}{n!} W_{n}^{Y} \left(\frac{x}{1-x}; r\right),
\end{equation}
where  $W_{n}^{Y} \left(\frac{x}{1-x}; r\right)$ are the $r$th order probabilistic Fubini polynomials coincides with (\ref{ar11}).\\
For a particular choice $f(x)=x^n$, from  (\ref{ni1}), we have 
\begin{equation}\label{eq123}
\sum_{i=0}^{\infty} \binom{-r}{i} \mathbb{E}S_i^n x^i = \frac{1}{(1-x)^r} W_{n}^{Y} \left(-\frac{x}{1+x}; r\right).
\end{equation}
This is a probabilistic extension of the  formula (3.28) studied in \cite{boyadzhiev2005series}.\\
Now, in the following propositions, we prove some intresting identities for the probabilistic Fubini polynomials.
\begin{proposition}\label{pr1}
	For $k \in \mathbb{N}_0$, we have
	\begin{equation}\label{de11}
	\frac{1}{1-x}\sum_{m=k}^{\infty} \binom{m}{k} \frac{u^m}{m!}v^{m-k}  W_m^Y\left(\frac{x}{1-x}\right) = \frac{u^k}{k!} D_{uv}^k \left(\frac{1}{1-x\mathbb{E}e^{uvY}}\right),
	\end{equation}
	provided $D_{t}^k \left(\displaystyle \frac{1}{1-x\mathbb{E}e^{tY}}\right)$ exists for $Y \in \mathcal{G},$ where $D_{t}^k$ is $k$th order differential operator with respect to $t$.	
\end{proposition}
\begin{proof}
	We start with left hand side of (\ref{de11})  and with the help of (\ref{eq1}), we get
	\begin{align*}
	\sum_{m=k}^{\infty} \binom{m}{k} \frac{u^m}{m!}v^{m-k}  W_m^Y\left(\frac{x}{1-x}\right) &= (1-x) \sum_{m=k}^{\infty} \binom{m}{k} \frac{u^m}{m!}v^{m-k} \sum_{n=0}^{\infty}\mathbb{E}S_n^m x^n\\
	&= (1-x)\mathbb{E}\left[\sum_{n=0}^{\infty} x^n 	\sum_{m=k}^{\infty} \binom{m}{k} \frac{u^m}{m!}v^{m-k} S_n^m\right]\\
	&= (1-x)\frac{u^k}{k!} \mathbb{E}\left[\sum_{n=0}^{\infty}x^n  S_n^k	e^{uvS_n}\right].
	\end{align*}
	Using $D_t^k\mathbb{E}e^{tS_n} = \mathbb{E}S_n^ke^{tS_n},$ we get
	\begin{align*}
	\frac{1}{1-x}\sum_{m=k}^{\infty} \binom{m}{k} \frac{u^m}{m!}v^{m-k}  W_m^Y\left(\frac{x}{1-x}\right) 
	&= \frac{u^k}{k!} \sum_{n=0}^{\infty}x^n D_{uv}^k\mathbb{E}e^{uvS_n}\\
	&= \frac{u^k}{k!} D_{uv}^k\left(\sum_{n=0}^{\infty}x^n \left(\mathbb{E}e^{uvY}\right)^n \right)\\
	&=  \frac{u^k}{k!} D_{uv}^k \left(\frac{1}{1-x\mathbb{E}e^{uvY}}\right).
	\end{align*}
	Hence, the identity is proved.
\end{proof}
\begin{remark}
	The Proposition \ref{pr1} may be viewed as a probabilistic extension to the identity proved in \cite{meoli2023some}.
\end{remark}

\begin{proposition} For $k \in \mathbb{N}$, we have
	\begin{equation}
	\frac{1}{1-x} \sum_{m = 1}^{\infty} \mathbb{E}S_m^k \frac{y^m}{m!}  W_m^Y\left(\frac{x}{1-x}\right) = \sum_{j=0}^{\infty} x^j B_k^Y \left(yD_y\right)\left(\mathbb{E}e^{yY}\right)^j,
	\end{equation}
	provided $D_y^i\left(\mathbb{E}e^{yY}\right)^j$ exists for $Y \in \mathcal{G}.$
\end{proposition}
\begin{proof}
	Using (\ref{eq1}), we get
	\begin{align*}
	\frac{1}{1-x} \sum_{m = 1}^{\infty} \mathbb{E}S_m^k \frac{y^m}{m!}  W_m^Y\left(\frac{x}{1-x}\right) 
	&= \sum_{m = 1}^{\infty} \mathbb{E}S_m^k \frac{y^m}{m!} \sum_{j = 0} x^j \mathbb{E}S_j^m\\
	&= \sum_{j=0}^{\infty} x^j \mathbb{E}\left[\sum_{m = 1}^{\infty} \left(\mathbb{E}S_m^k \right) \frac{(yS_j)^m}{m!}\right]\\
	&= \sum_{j=0}^{\infty} x^j \mathbb{E}\left[e^{yS_j} B_k^Y(yS_j)\right], \;\;\; (\text{using } (\ref{0009})),\\
	&= \sum_{j=0}^{\infty} x^j \sum_{i = 0}^{k}S_Y(k,i)y^i D_y^i\left(\mathbb{E}e^{yS_j}\right).
	\end{align*}
	With the help of (\ref{bp11}), we get the proposition. 
\end{proof}

\end{example}
\section{Probabilistic Fubini Polynomials and Some Probability Distributions} \noindent For  different choices of probability distribution of the rv $Y$, we  obtain the different representations of the probabilistic Fubini polynomials and the numbers. These representations may be in terms of the Stirling numbers of the second kind,  polylogarithm functions, and   the Apostol-Euler polynomials.  Some spacial choices of the rvs $Y \in \mathcal{G}$, we have the following examples.

\begin{example}
Let $Y$ be a Poisson random variate with $pmf$  and the moment generating function ($mgf$) %
$\mathbb{P}\{Y = k\} = e^{-\lambda}\frac{\lambda^k}{k!},$ 
and
$\mathbb{E}\left(e^{tY}\right) = e^{\lambda(e^t -1)},\;\lambda>0$,
respectively.\\
Substituting the $mgf$ into (\ref{pp1}) and using (\ref{0003}), we obtain

\begin{align*}
\sum_{n=0}^{\infty} W_n^Y(x) \frac{t^n}{n!} 
= \frac{1}{1-x\left(e^{\lambda (e^t -1)}-1\right)}
&= \sum_{i=0}^{\infty}\lambda^i W_i(x) \frac{(e^t -1)^i}{i!}
= \sum_{n=0}^{\infty} \left(\sum_{i=0}^{n} \lambda^i W_i(x)S(n,i)\right)\frac{t^n}{n!}.
\end{align*}
Finally, comparing the coefficients of powers of $t$, we get a convolution result of the form
\begin{equation*}\label{55}
W_n^Y(x) = \sum_{i=0}^{n} \lambda^i W_i(x)S(n,i).
\end{equation*}
Let $Y_1, Y_2, \dots, Y_i$ be  IID copies  of Poisson rv with mean $\lambda$. Then, $S_i = Y_1 + Y_2 +\cdots + Y_i\sim \text{Poisson}(i\lambda)$ for $i = 1,2,\dots$. Using (\ref{00001}) and (\ref{eq1}), we obtain a relationship between the Bell polynomials and the probabilistic Fubini polynomials as  
\begin{equation*}\label{551}
\sum_{i=0}^{\infty} B_n(i \lambda) x^i = \frac{1}{1-x} W_n^Y \left(\frac{x}{1-x}\right) = \sum_{k=0}^{\infty}x^k \sum_{j=0}^{n \wedge k}\binom{k}{j} j! S_Y(n,j).
\end{equation*}
\end{example}

\begin{example}\label{exp42}

Let $Y$ be a geometric rv different from $G_p$ with the $pmf$
\begin{equation*}
\mathbb{P}\{Y = k\} = rs^{k-1}, \;\;\; k = 1,2,\dots,
\end{equation*}
where $s =1-r,\;\;0<r\leq 1.$

One can  verify the following interconnection between the polylogarithms and the geometric  variate
\begin{equation*}
Li_{-n}(s) = \frac{s}{r} \mathbb{E}Y^n,
\end{equation*}
where $Li_z(y) = \sum_{i = 1}^{\infty} \frac{y^j}{i^z}$ with $y,z \in \mathbb{C}$ and $|y| <1$.

We define the $k$th multinomial convolution of the polylogarithm function as

\begin{equation*}
Li_{-n}^{*k}(s) = \sum_{n_1 +n_2 + \cdots +n_k =n}^{}\frac{n!}{n_1! n_2! \cdots n_k!}Li_{-n_1}(s)Li_{-n_2}(s)\cdots Li_{-n_k}(s),~~ k \in \mathbb{N},
\end{equation*}
with $Li_0^{*0}(s) = 1$.

Let $\langle  Y_i \rangle _ {i\geq 0}$ be the sequence of independent copies of the geometric rv $Y$. Then, we obtain

\begin{equation}\label{789}
Li_{-n}^{*k}(s) = \frac{s^k}{r^k}  \mathbb{E}S_k^n,   \;\;\; \forall \; k \in \mathbb{N}.
\end{equation}

\begin{theorem} Let $Y$ be the geometric rv as considered in Example \ref{exp42}. Then
\begin{equation}\label{1012}
\mathbb{E}\left(\left(\frac{r}{s}\right)^{G_{p}-1}Li_{-n}^{*(G_{p}-1)}(s)\right) = \lim_{m \rightarrow \infty} \sum_{k=0}^{m} r^{k+1} Li_{-n}^{*k}(s) = W_n^Y(x),
\end{equation}
where $G_p$ is geometric variate with parameter $p =\eta(x)= 1/(1+x)$ as defined in (\ref{grv1}).\\ Also, for $p = (1+x)/(1+2x),$ we have
\begin{equation}\label{eq2}
\lim_{m \rightarrow \infty} \sum_{k=0}^{m} \left(\frac{r}{s}\right)^k Li_{-n}^{*k}(s) x^k = (1+x) W_n^Y(x),
\end{equation}
provided above limits exist.
\end{theorem}
\begin{proof}
On multiplying with $pq^k$ both sides of (\ref{789}) and using (\ref{grv1}), we get
\begin{align}\label{eq3}
\frac{r^k}{s^k} Li_{-n}^{*k}(s) pq^k &=   \mathbb{E}S_k^n pq^k, \;\; \forall k \in \mathbb{N} \nonumber \\
\lim_{m \rightarrow \infty} \sum_{k=0}^{m} \frac{r^k}{s^k} Li_{-n}^{*k}(s) pq^k &=  \lim_{m \rightarrow \infty} \sum_{k=0}^{m} \mathbb{E}S_k^n pq^k,
\end{align}
where $q=1-\eta(x)$. For $p = \eta(x),$ the right hand side quantity of (\ref{eq3}) converges to the $n$th order moment of geometric rv provided the limit exists. Hence, using (\ref{ad}), we get the desired result.
The result (\ref{eq2}) is the consequence of identity (\ref{eq1}) with the help of (\ref{1012}).
\end{proof}

\end{example}

\begin{example}
Let $Y$ be a Bernoulli rv  with
$\mathbb{E}(e^{tY})-1 = p(e^t-1),\;0<p\leq 1$.
Clearly, from (\ref{pp1}), we get
\begin{equation*}
W_n^Y(x) = W_n(px).
\end{equation*}
On the other hand, the Apostol-Euler polynomials $E(c;x)$ are defined by the exponential generating function of the form (see \cite{adell2019closed})
\begin{equation}\label{aep}
\frac{e^{xt}}{1+c(e^t - 1)} = \sum_{n=0}^{\infty} E(c;x) \frac{t^n}{n!},
\end{equation}
where $t \in \mathbb{R}$ and $  c \in [0,1].$\\
For $-1 \leq cp \leq 0,$ consider  (\ref{aep}) and with the help of (\ref{pp1}), we get

\begin{align*}
\sum_{n=0}^{\infty} E(-cp;x) \frac{t^n}{n!}	&
= \left(\sum_{k=0}^{\infty} W_k^Y(c) \frac{t^k}{k!}\right) \left(\sum_{i=0}^{\infty} x^i \frac{t^i}{i!} \right)
= \sum_{n=0}^{\infty} \left(  \sum_{k=0}^{n} \binom{n}{k} W_k^Y(c)x^{n-k}\right) \frac{t^n}{n!}.
\end{align*}
We establish an interconnection between the Apostol-Euler polynomials and the probabilistic Fubini polynomials by comparing the coefficients of $t$. It is given by
\begin{equation}\label{aep1}
E(-cp;x) = \sum_{k=0}^{n} \binom{n}{k} W_k^Y(c)x^{n-k}.
\end{equation}
\end{example}

\begin{remark}
For $x=0$, (\ref{aep1}) gives
\begin{equation*}
E(-cp) = W_n^Y(c),
\end{equation*}
where $E(\cdot)$ are the Apostol-Euler numbers.
\end{remark}
\noindent For any $\alpha \in \mathbb{R}$, we have established an interconnection between higher order probabilistic Fubini polynomials and the generalized Apostol-Euler polynomials  of the following form
\begin{equation*}
E(\alpha,-cp;x) = \sum_{k=0}^{n} \binom{n}{k} W_k^{Y}(c;\alpha)x^{n-k},
\end{equation*} 
where $E(\alpha,-cp;x)$ are the generalized Apostol-Euler polynomials defined in \cite{adell2019closed}.

\section{Probabilistic Fubini Numbers and Its Determinant Expressions}
Recently, the determinant expressions of several polynomials and numbers are obtained in the literature. Komatsu \cite{komatsu2017hypergeometric} and Glaisher \cite{glaisher1876expressions} studied the determinant expressions for the Cauchy polynomials, Bernoulli numbers and the Euler numbers.\\
The following theorem provides a determinant expression for a sequence of the real numbers.
\begin{theorem}\label{eqtm}
(Komatsu \cite{KomatsuRamrez}) Let  $\langle f(n) \rangle_{n \in \mathbb{N}}$ be a sequence with $f(0)=1$ and let $w(k)$ be an arbitrary function independent of $n$. Then
\begin{equation*}
f(n) = \begin{vmatrix}
w(1) & 1 & 0 & 0& \cdots & 0 \\ 
\\
w(2) & w(1) & 1 &0&  \cdots & 0\\
\\
w(3) & w(2) & w(1) & 1  & \cdots & 0 \\
\\
\vdots & \vdots & \vdots & \vdots& \ddots & \vdots \\
\\
w(n-1) & w(n-2) &w(n-3)&w(n-4)& \cdots & 1\\
\\
w(n) & w(n-1)& w(n-2) &w(n-3) & \cdots & w(1)
\end{vmatrix}
\end{equation*}
if and only if
\begin{equation}\label{eqq2}
f(n) = \sum_{k=1}^{n} (-1)^{k-1} w(k)f(n-k) \text{ with } n \geq 1.
\end{equation}
Also, function $w(k)$ is expressed as
\begin{equation}\label{det1}
w(k) = \begin{vmatrix}
f(1) & \;\;\;\;1 & 0 & 0& \cdots & 0 \\ 
\\
f(2) & \;\;f(1) & 1 & 0& \cdots & 0\\
\\
f(3) & f(2) & f(1) & 1&\cdots & 0\\
\\
\vdots & \vdots & \vdots & \vdots & \ddots&\vdots \\
\\
f(k-1)  & f(k-2)&f(k-3)&f(k-4)& \cdots&1\\
\\
f(k)  & f(k-1) &f(k-2)&f(k-3) & \cdots& f(1)\\
\end{vmatrix}.
\end{equation}
\end{theorem} 

The following lemma (see \cite{KomatsuRamrez, komatsu2017hypergeometric, merca2013note}) will be used to obtain the explicit expression for sequence of the probabilistic Fubini numbers.
\begin{lemma}\label{lemma01}
Let $A$ be a square matrix of order $(k+1)$ defined by
\begin{equation*}
A = 	\begin{bmatrix}
1 &0 & \cdots & 0 \\
\\
f(1) & 1 & \cdots & 0 \\ 
\\
\vdots & \vdots & \ddots & \vdots \\
\\
f(k)  &f(k-1)& \cdots & 1\\
\end{bmatrix}.
\end{equation*}
Also, inverse of $A$ is given by
\begin{equation*}
A^{-1} = 	\begin{bmatrix}
1  &0 & \cdots & 0\\
\\
w(1) & 1  & \cdots & 0 \\ 
\\
\vdots & \vdots & \ddots&\vdots \\
\\
w(k)  &w(k-1)& \cdots  & 1\\
\end{bmatrix}.
\end{equation*}
Using Trudi's formula (see \cite{KomatsuRamrez, merca2013note}), the combinatorial expression of sequence $f(n)$ is obtained. It has the following combinatorial form
\begin{equation*}
f(n) = \sum_{l_1 + 2l_2 + \cdots + nl_n = n}^{}\binom{l_1 + \cdots + l_n}{l_1,\dots, l_n}(-1)^{n-l_1 - \cdots -l_n} w(1)^{l_1}w(2)^{l_2} \cdots w(n)^{l_n},
\end{equation*} 
where $\binom{l_1 + \cdots + l_n}{l_1,\dots, l_n}$ are multinomial coefficients and $l_i$'s stand for the numbers of blocks with $i$ elements while partitioning a set with $n$ elements.
\end{lemma}
In the next result, we obtain a determinant expression to the probabilistic Fubini numbers  and present a combinatorial sum formula for the probabilistic Fubini numbers. 

\begin{theorem}
For $n \geq 1$, we have
\begin{equation}\label{det1}
W_n^Y = n!	\begin{vmatrix}
\frac{\mathbb{E}Y}{1!} & 1  &0 & 0&\cdots & 0 \\ 
\\
-\frac{\mathbb{E}Y^2}{2!} & \frac{\mathbb{E}Y}{1!} & 1 &0 & \cdots & 0\\
\\
\frac{\mathbb{E}Y^3}{3!}&-\frac{\mathbb{E}Y^2}{2!} & \frac{\mathbb{E}Y}{1!} & 1 &\cdots & 0\\
\\
\vdots & \vdots & \vdots & \vdots&\ddots &\vdots \\
\\
(-1)^{n-2}\frac{\mathbb{E}Y^{n-1}}{(n-1)!} & (-1)^{n-3}\frac{\mathbb{E}Y^{n-2}}{(n-2)!} & (-1)^{n-4}\frac{\mathbb{E}Y^{n-3}}{(n-3)!} & (-1)^{n-5}\frac{\mathbb{E}Y^{n-4}}{(n-4)!}& \cdots & 1\\
\\
(-1)^{n-1}\frac{\mathbb{E}Y^n}{n!} & (-1)^{n-2}\frac{\mathbb{E}Y^{n-1}}{(n-1)!} & (-1)^{n-3}\frac{\mathbb{E}Y^{n-2}}{(n-2)!} &(-1)^{n-4}\frac{\mathbb{E}Y^{n-3}}{(n-3)!} & \cdots & \frac{\mathbb{E}Y}{1!}
\end{vmatrix}.
\end{equation}
Moreover, it has a explicit combinatorial expression of the form
\begin{equation}\label{mws1}
W_n^Y = n!\sum_{l_1 + 2l_2 + \cdots + nl_n = n}^{}\binom{l_1 + \cdots + l_n}{l_1,\dots, l_n} \left(\frac{\mathbb{E}Y}{1!}\right)^{l_1}\left(\frac{\mathbb{E}Y^2}{2!}\right)^{l_2} \cdots \left(\frac{\mathbb{E}Y^n}{n!}\right)^{l_n}.
\end{equation}
\end{theorem}

\begin{proof}
Simplifying the recurrence relation obtained in Proposition \ref{p3} for $x=1$, we get
\begin{equation}\label{eqq1}
\frac{W_n^Y}{n!} = \sum_{k=1}^{n} \frac{\mathbb{E}Y^k}{k!} \frac{W_{n-k}^Y}{(n-k)!}.
\end{equation}
Observe that (\ref{eqq1}) has a similar expression as (\ref{eqq2}) with 
\begin{equation}\label{wk}
f(n) = \frac{W_n^Y}{n!} \text{ and } w(k) = (-1)^{k-1} \frac{\mathbb{E}Y^k}{k!}.
\end{equation}
Using Theorem \ref{eqtm}, the required determinant expression (\ref{det1}) can be obtained.

From (\ref{wk}), we also have 
\begin{equation*}
\begin{bmatrix}
1 &  0 & 0 &0& \cdots & 0 \\
\\
\frac{\mathbb{E}Y}{1!} & 1 &0  &0& \cdots &0\\ 
\\
-\frac{\mathbb{E}Y^2}{2!} & \frac{\mathbb{E}Y}{1!} & 1& 0 & \cdots & 0\\
\\
\vdots & \vdots & \vdots & \vdots & \ddots & \vdots \\
\\
w(n-1) & w(n-2) & w(n-3) & w(n-4) & \cdots & 0\\
\\
w(n) & w(n-1) & w(n-2) & w(n-3) & \cdots & 1
\end{bmatrix}^{-1} = 		\begin{bmatrix}
1  &  0 &0 & 0& \cdots & 0\\
\\
\frac{W_1^Y}{1!} & 1  &0  &0& \cdots & 0\\ 
\\
\frac{W_2^Y}{2!} & \frac{W_1^Y}{1!} & 1 & 0 & \cdots & 0 \\
\\
\vdots & \vdots & \vdots & \vdots & \ddots & \vdots \\
\\
\frac{W_{n-1}^Y}{(n-1)!} & \frac{W_{n-2}^Y}{(n-2)!} &  \frac{W_{n-3}^Y}{(n-3)!} &  \frac{W_{n-4}^Y}{(n-4)!} & \cdots&0\\
\\
\frac{W_n^Y}{n!} & 	\frac{W_{n-1}^Y}{(n-1)!} &\frac{W_{n-2}^Y}{(n-2)!} & \frac{W_{n-3}^Y}{(n-3)!} & \cdots & 1
\end{bmatrix}.
\end{equation*}
Hence, using Lemma \ref{lemma01} and Trudi's formula (see \cite{KomatsuRamrez, merca2013note}), we get the required combinatorial interpretation of the probabilistic Fubini numbers.
\end{proof}
For degenerate rv $Y$ at $1$, we get the combinatorial interpretation of the Fubini numbers as (see \cite{KomatsuRamrez})
\begin{equation*}
W_n = n!\sum_{l_1 + 2l_2 + \cdots + nl_n = n}^{}\binom{l_1 + \cdots + l_n}{l_1,\dots, l_n}\left(\frac{1}{1!}\right)^{l_1}\left(\frac{1}{2!}\right)^{l_2} \cdots \left(\frac{1}{n!}\right)^{l_n}.
\end{equation*}
\begin{example}
Suppose  $Y$ follows standard exponential distribution. Then, using (\ref{det1}) and (\ref{wk}), we get
\begin{equation*}
(-1)^{n-1} = 	\begin{vmatrix}
\frac{W_1^Y}{1!} & 1&0 &0&\cdots&0\\ 
\\
\frac{W_2^Y}{2!} & \frac{W_1^Y}{1!} & 1&0&\cdots& 0 \\
\\
\frac{W_3^Y}{3!}&\frac{W_2^Y}{2!} & \frac{W_1^Y}{1!} & 1&\cdots& 0 \\
\\
\vdots & \vdots & \vdots & \vdots &\ddots &\vdots\\
\\
\frac{W_{n-1}^Y}{(n-1)!} & \frac{W_{n-2}^Y}{(n-2)!} & \frac{W_{n-3}^Y}{(n-3)!} & \frac{W_{n-4}^Y}{(n-4)!}&\cdots & 1\\
\\
\frac{W_n^Y}{n!} & 	\frac{W_{n-1}^Y}{(n-1)!} & \frac{W_{n-2}^Y}{(n-2)!} &\frac{W_{n-3}^Y}{(n-3)!}& \cdots & \frac{W_1^Y}{1!}
\end{vmatrix}.
\end{equation*}
Also, from (\ref{mws1}), we have weighted sum of multinomial coefficients as
\begin{equation*}
W_n^Y = n!\sum_{l_1 + 2l_2 + \cdots + nl_n = n}^{}\binom{l_1 + \cdots + l_n}{l_1,\dots, l_n}.
\end{equation*}

\end{example}

\bibliographystyle{acm}
\bibliography{bbbb11}

\end{document}